\documentclass[11pt]{article}
\usepackage{latexsym}
\usepackage{amsmath}
\usepackage{graphicx} 
\usepackage{amsfonts} 
\usepackage{amssymb}
\newtheorem{thm}{Theorem}[section]
\newtheorem{pro}[thm]{Proposition}
\newtheorem{lem}[thm]{Lemma}
\newtheorem{defn}[thm]{Definition}
\newtheorem{con}[thm]{Conjecture}
\newtheorem{cor}[thm]{Corollary}
\newtheorem{rmk}[thm]{Remark}

\newcommand{\mc}{\left[ \begin{array}{c} \widetilde{m} \\
\left\lfloor \frac{|I_{1}|+1}{2} \right\rfloor , \ldots , 
\left\lfloor \frac{|I_{s}|+1}{2} \right\rfloor \end{array} \right]_{x^2}}
\newcommand{\mcb}{\left[ \begin{array}{c} \widetilde{m} \\
\left\lfloor \frac{|J_{1}|+1}{2} \right\rfloor , \ldots , 
\left\lfloor \frac{|J_{s}|+1}{2} \right\rfloor \end{array} \right]_{x^2}}
\newcommand{\p}{\noindent}
\newcommand{\qa}{S_n^{I}}
\newcommand{\qb}{B_n^{J}}
\newcommand{\mt}{\widetilde{m}}

\newcommand{\sumb}{\sum_{\sigma \in {\qb}} (-1)^{\ell(\sigma)} x^{L(\sigma)}}
\newcommand{\N}{\mathbb N}

\begin{document}
\begin{center}
{\Large \bf Proof of a conjecture of Klopsch-Voll on \\
 Weyl groups of type $A$ 
}
 \vspace{0.8cm} \\
 Francesco Brenti \\
Dipartimento di Matematica \\
Universit\`{a} di Roma ``Tor Vergata''\\
Via della Ricerca Scientifica, 1 \\
00133 Roma, Italy \\
{\em brenti@mat.uniroma2.it } \\
 \vspace{0.5cm}
 Angela Carnevale \footnote{Current address: Fakultat f\"ur Mathematik,
Universit\"at Bielefeld,
D-33501 Bielefeld, Germany, {\em acarneva1@math.uni-bielefeld.de}.} \\
Dipartimento di Matematica \\
Universit\`{a} di Roma ``Tor Vergata''\\
Via della Ricerca Scientifica, 1 \\
00133 Roma, Italy \\
{\em carneval@mat.uniroma2.it } \\
\end{center}
\vspace{1cm}
\begin{abstract}
We prove a conjecture of Klopsch-Voll on the signed generating function of a 
new statistic on the quotients of the symmetric groups. As a consequence of our
results we also prove a conjecture of Stasinski-Voll in type $B$.
\end{abstract}

\section{Introduction}

In \cite{KV} a new statistic on Weyl groups of type $A$ was introduced, which combines combinatorial
and parity conditions, in connection with formed spaces. In the same paper
 the authors conjecture a relationship between the signed (by length)
generating function of this new statistic over the quotients of the symmetric groups and the enumeration 
of partial flags in a non-degenerate quadratic space  (see \S 1.2.1 and Conjecture C of \cite{KV}, for details).

The purpose of this work is to prove this conjecture. As a consequence of our results we also prove a conjecture in
\cite{VS2} (see Conjecture 1.6) which relates the generating function of an analogous statistic over the quotients
of the Weyl groups of type $B$  and the local factors of the representation zeta function of certain groups (see \S 1.3
and Theorem C of \cite{VS2}, for details). Our proofs are combinatorial.

The organization of the paper is as follows. In the next section we collect some notation, definitions, and results
that are used in the sequel. In \S 3 we prove some preliminary lemmas that are used in \S 4 in the proof of
our main result. In particular, we prove that certain operations on a quotient of the symmetric group do not 
change the relevant generating function. In \S 4, using these results, we prove our main result, namely
Conjecture C of \cite{KV}. In \S 5, as a consequence of our main result, we prove Conjecture 1.6 of \cite{VS2}.

After we submitted this paper we learned that Conjecture 1.6 of 
\cite{VS2} has also been proved, independently and in a different way, by A. Landesman in \cite{Lan}. Landesman's proof shows that
both sides of the identity given by Conjecture 1.6 of \cite{VS2} satisfy the same recursion.

\section{Preliminaries}

In this section we introduce some notation, definitions, and results that are used in the sequel. 

\p
We let $\mathbb P:=\{1,\,2,\ldots\}$ be the set of positive integers and $\N:= \mathbb{P} \cup \{0\}$. For all $m,\,n \in \mathbb{Z}$, $m\leq n$ we let $[m,n]:= \{m,\,m+1,\ldots,\,n\}$ and $[n]:=[1,\,n]$. Given a set $I$ we denote by $|I|$ its cardinality.

\p
We follow $\cite{BB}$ for notation and terminology about Coxeter groups. 

The symmetric group $S_n$ is the group of permutations of the set $[n]$. For $\sigma \in S_n$ we use both the one-line notation $\sigma=[\sigma(1),\,\ldots,\,\sigma(n)]$ and the disjoint cycle notation.
We let $s_1,\ldots,\,s_{n-1}$ denote the standard generators of $S_n$, $s_i=(i,\,i+1)$.

The hyperoctahedral group $B_n$ is the group of signed permutations, or permutations $\sigma$ of the set $[-n,n]$ such that $\sigma(j)=-\sigma(-j)$. 
Given such a $\sigma$ we write $\sigma = [a_1, \ldots ,a_n]$ to mean that $a_j = \sigma(j)$ for $j=1, \ldots ,n$.
The Coxeter generating set of $B_n$ is $S=\{s_0,\,s_1,\,\ldots,\,s_{n-1}\}$, where $s_0=[-1,\,2,\,3,\ldots,\,n]$ and $s_1,\ldots,\,s_{n-1}$ are as above. 

For $(W,S)$ a Coxeter system we let $\ell$ be the Coxeter length and for 
$I\subseteq S$ we define the quotients:
\[
W^{I} := \{w\in W \;:\; D(w)\subseteq S\setminus I\},
\]
and 
\[
^{I}W := \{w\in W \;:\; D_L(w)\subseteq S\setminus I\},
\]
where $D(w)=\{s\in S \;:\; \ell(ws)<\ell(w)\}$, and 
$D_L(w)=\{s\in S \;:\; \ell(sw)<\ell(w)\}$, and the parabolic subgroup $W_I$
to be the subgroup generated by $I$. For subsets $X\subseteq W$ we let $X^{I}:= X \cap W^{I}$. The following result is well known (see, e.g., \cite[Proposition 2.4.4]{BB}).
\begin{pro}
Let $(W,S)$ be a Coxeter system, $J \subseteq S$, and $w \in W$. Then there exist
unique elements $w^J \in W^J$ and $w_J \in W_J$ (resp., $^Jw \in ^JW$ and 
$_Jw \in W_J$) such that $w= w^J w_J$ (resp., $_Jw ^Jw$). Furthermore
$\ell(w)= \ell(w^J)+\ell(w_J)$ (resp., $\ell(_Jw)+\ell(^Jw)$).
\end{pro}

It is well known that $S_n$ and $B_n$, with respect to the above generating 
sets, are Coxeter systems and that the following results hold (see, e.g., 
\cite[Propositions 1.5.2, 1.5.3, and \S 8.1]{BB}).

\begin{pro}
Let $\sigma \in S_n$. Then
\[
\ell(\sigma)=| \{ (i,j) \in [n]^2 : i<j , \sigma(i) > \sigma(j) \} |
\]
and
\[
D(\sigma) = \{ s_i :  \sigma(i) > \sigma(i+1) \}.
\]
\end{pro}

\begin{pro}
Let $\sigma \in B_n$. Then
\[
\ell(\sigma)= \frac{1}{2} | \{ (i,j) \in [-n,n]^2 : i<j , \sigma(i) > \sigma(j) \} |
\]
and
\[
D(\sigma) = \{ s_i : i \in [0,n-1] , \sigma(i) > \sigma(i+1) \}.
\]
\end{pro}

The following statistic was first defined in \cite{KV}, see Definition 5.1.
\begin{defn}
Let $n\in {\mathbb P}$. The statistic $L_{A}:S_n \rightarrow \N$ is defined as follows. For $\sigma \in S_n$ 
\[
L_{A}(\sigma):= \sum_{I \subseteq [n-1]} (-1)^{|I|} 2^{n-2-|I|} \ell(^I \sigma).
\]
\end{defn}

The following result is proved in \cite[Lemma 5.2]{KV}.

\begin{lem}
Let $n\in {\mathbb P}$ and $\sigma \in S_n$. Then  
\[
L_{A}(\sigma) = |\{(i,j) \in [n]^2 \; : \; i<j,\,\sigma(i)>\sigma(j),\,i\not\equiv j \pmod{2}\}|.
\]
\end{lem}
For example let $n=5$, $\sigma=[4,2,1,5,3]$. Then $L_A (\sigma)=|\{(1,2),\,(2,3),\,(4,5)\}|=3$.

Following \cite{KV} and \cite{VS}, we define {\em chessboard} elements, both in $S_n$ and $B_n$, as follows. 

Let $n\in {\mathbb P}$ and $W$ be $S_n$ or $B_n$. Set:\begin{align*}C_{n,\, +} &:=  \{w\in W \;:\; i+ w(i)\equiv 0 \pmod{2},\;  i=1,\ldots,\,n \} \\C_{n,\,-} &:=  \{w\in W \;:\; i+ w(i)\equiv 1 \pmod{2},\;  i=1,\ldots,\,n \} \\C_n       &:=  C_{n,\,+} \cup C_{n,\,-}.\end{align*}\p For $n=2m+1$ clearly $C_{n,\,-}=\emptyset$ so $C_n = C_{n,\,+}$. Note that the chessboard elements $C_n$ form a subgroup of $W$ and the {\em even} chessboards elements $C_{n,\,+}$ form a subgroup of $C_n$. Thus on chessboard elements one can define, besides the usual sign charachter (restriction of the sign character over $W$), the linear character $\chi : C_n \longrightarrow \{\pm 1\}$, whose kernel is the subgroup of even chessboard elements, $\ker \chi=C_{n,\,+}$.Clearly $\chi$ is trivial over $C_n$ for $n$ odd.

For a real number $x$ we denote by $\left\lfloor x \right\rfloor$ the greatest integer less than or equal to $x$ and by $\left\lceil x \right\rceil$ the smallest integer greater than or equal to $x$.

\p
Finally, for $n_1,\ldots,\,n_k \in \N$ such that $\sum_{i=1}^k n_i =n$ we denote by \footnotesize $\left[ \begin{array}{c} n \\ n_1,\ldots,n_k\end{array} \right]_{q}$ \normalsize the {\em $q-$multinomial coefficient} 
\[\left[ \begin{array}{c} n \\
n_1,\ldots,n_k \end{array} \right]_{q}:= \frac{[n]_{q}!}{[n_1]_q !\cdot \ldots\cdot [n_k]_q !},\]
where 
\[[n]_q := \frac{1-q^n}{1-q},\qquad\qquad [n]_q ! := \prod_{i=1}^{n} [i]_q \qquad\qquad\mbox{and}\qquad [0]_q!:= 1.
\]

Given $J \subseteq [0,n-1]$ there are unique integers $a_1 < \cdots < a_s$ and  $b_1 < \cdots < b_s$ such that $J = [a_1,b_1] \cup \cdots \cup [a_s,b_s]$ and  $a_{i+1} - b_{i} >1$ for $i=1, \ldots , s-1$. We call the intervals  $[a_1,b_1], \ldots , [a_s,b_s]$ the {\em connected components} of $J$.
The following conjecture was made in \cite[Conjecture C]{KV}.

\begin{con}
\label{A}
Let $n\in {\mathbb P}$, $m:=\left\lfloor{\frac{n}{2}}\right\rfloor$ and $I \subseteq [n-1]$. Then
\begin{align}\small\sum_{\sigma \in C^{I}_{2m+1}} (-1)^{\ell (\sigma )} \chi(\sigma) x^{L(\sigma )} = &\mc \displaystyle\prod\limits_{\tiny{k=\mt + 1}}^{\tiny m}{(1-x^{2k})} \\\sum_{\sigma \in C^{I}_{2m}} (-1)^{\ell (\sigma )} \chi(\sigma) x^{L(\sigma )} = &\left\{ \begin{array}{ll}\mc , & \mbox{if $m=\mt$,} \\ (1+x^{m}) \mc \displaystyle\prod\limits_{\tiny{k=\mt + 1}}^{\tiny {m-1}}{(1-x^{2k})}, & \mbox{otherwise,}\end{array} \right.  
\end{align}\normalsize 
where $I_{1},\ldots , I_{s}$ are the connected components of $I$ and $\tilde{m} := \sum_{k=1}^{s} \lfloor \frac{|I_{k}|+1}{2} \rfloor$.
\end{con}

\begin{rmk}
One can check, using Proposition 4.5 in \cite{KV}, that Conjecture \ref{A} is  indeed equivalent to Conjecture C of \cite{KV}.  
\end{rmk}

Conjecture \ref{A} is known to be true if $|I| \geq n-2$ (see \cite[p. 4433]{KV}).
The purpose of this work is to prove Conjecture \ref{A} in full generality.
As a consequence of our results we also prove a conjecture in \cite{VS2}, which
we now describe.

The following statistic was introduced in \cite{VS} and \cite{VS2}, and is a natural
analogue of the statistic $L_A$ introduced above, for Weyl groups of type $B$.

\begin{defn}
\label{defLB}
Let $n\in {\mathbb P}$. The statistic $L_{B}:B_n \rightarrow \N$ is defined as follows. For $\sigma \in B_n$ 
\[
L_{B}(\sigma):= \frac{1}{2} |\{(i,j)\in [-n,\,n]^2 \, : \, i<j,\,\sigma(i)>\sigma(j),\,i\not\equiv j \pmod{2}\}|.\]
\end{defn}
For example, if $n=4$ and $\tau=[-2,4,3,-1]$ then $L_B (\tau)= \frac{1}{2}|\{(-4,-3),\,(-4,1),\,(-3,-2),$ $\,(-1,0),\,(-1,4),\,(0,1),\,(2,3),\,(3,4)\}|=4$. 

We call these statistics $L_A$ and $L_B$ the {\em odd length} of the symmetric and hyperoctahedral groups, respectively.
Clearly $L_B(\mbox{id})=0$, while $L_B(s_i)=1$, for $i=0,\,1,\ldots,n-1$. 
Note that if $\sigma \in S_n \subset B_n $ then $ L_B (\sigma)=L_A(\sigma)$, so in the following we omit the subscript and write just $L$ for both statistics.


The following conjecture was made in \cite[Conjecture 1.6]{VS2}.

\begin{con}
\label{B}
Let $n\in {\mathbb P}$ and $J\subseteq [0,\,n-1]$. Then
\small
\begin{equation*}\sumb=\frac{\prod_{j=a+1}^{n}(1-x^i)}{\prod_{i=1}^{\mt}(1-x^{2i})}  \mcb 
\end{equation*}
\normalsize
\noindent 
where $J_0$ is the (possibly empty) connected component containing $0$, $J_1,\ldots,\,J_s$ are the remaining connected components of $J$, $\mt:= \sum_{i=1}^{s} \left\lfloor \frac{|J_i|+1}{2}\right\rfloor$ and $a:= \min\{ [0, \,n] \setminus J\}$.
\end{con}

Conjecture \ref{B} is known to hold if $J=[n-1]$, if $J= \emptyset$, and if
$n \equiv 0 \pmod{2}$ and $[0,n-1] \setminus J \subseteq 2 {\mathbb N}$
(see \cite[Theorem 2]{VS}). In particular, the following holds
(see \cite[Proposition 9]{VS}).

\begin{pro}
\label{ascending}
Let $n\in {\mathbb P}$. Then 
\begin{equation*}
\sum_{\sigma \in {B_n^{[n-1]}}}{(-1)^{\ell(\sigma)} x^{L(\sigma)}}=
\prod_{j=1}^{\lceil \frac{n}{2} \rceil}{(1-x^{2j-1})}.
\end{equation*}
\end{pro}

We conclude with three results that are used in the sequel. The first one is proved in \cite{VS} (see Lemma 8).

\begin{lem}
\label{chessB}
Let $J \subseteq [0,n-1]$. Then
\[\sumb=\sum_{\sigma \in C^{J}_{n,\, +}}(-1)^{\ell (\sigma )}x^{L(\sigma )}.\]
\end{lem}

The next result follows easily from Corollary 20 and Proposition 22 of \cite{VS}.

\begin{pro}\label{fact2}
Let $n\in {\mathbb P}$ and $J\subseteq [n-1]$. If $n\equiv 1 \pmod 2$ or $n\equiv 0 \pmod 2$ and $[n-1]\setminus J \subseteq 2\N$   then
\[\sumb=\left(\sum_{\sigma \in B^{[n-1]}_{n}} (-1)^{\ell(\sigma)} x^{L(\sigma)}\right)\left(\sum_{\sigma \in S^{J}_{n}} (-1)^{\ell(\sigma)} x^{L(\sigma)}\right).\]
\end{pro}
The following result follows from the proof of Proposition 25 of \cite{VS}.
\begin{pro}
\label{2.8b}
 Let $n\in {\mathbb P}$ be even, and $J \subseteq [0,n-1]$ be such that $[0,n-1] \setminus
 J \subseteq 2 {\mathbb N}$.  Then
\[\sum_{\sigma \in B^{J {\scriptsize\setminus\{0\}}} _n } 
(-1)^{\ell (\sigma )} x^{L(\sigma )} = \left(\sumb\right) \left(\sum_{\sigma \in B^{[i-1]}_{i}} (-1)^{\ell(\sigma)} x^{L(\sigma)}\right), \]
where $i := min \{ [0,n] \setminus J \}$.
\end{pro}

\section{Shifting and compressing}

In this section we prove some results that are used in the next one in the proof of 
our main result. In particular, we prove that certain operations on a quotient of the symmetric group do not 
change the corresponding generating function.

The following result is the analogue, in type $A$, of Lemma \ref{chessB}.
\begin{lem}
\label{chessA}
Let $n \in  {\mathbb P}$ and $I \subseteq [n-1]$. If $n \equiv 1 \pmod{2}$ or $n \equiv 0 \pmod{2}$ and $[n-1] \setminus I 
\subseteq 2 \, {\mathbb N}$ then
\[ \sum_{\sigma \in S_{n}^{I}}(-1)^{\ell (\sigma )}x^{L(\sigma )} = \sum_{\sigma \in C_{n,+}^{I}}(-1)^{\ell (\sigma )}x^{L (\sigma )}. \]
\end{lem}
{\bf Proof:} Reasoning as in the proof of Lemma \ref{chessB} (i.e., of Lemma 8 in \cite{VS}) one can see that
\begin{equation}
\label{chessA0}
\sum_{\sigma \in S_{n}^{I}}(-1)^{\ell (\sigma )} \, x^{L(\sigma )}= \sum_{\sigma \in C_{n}^{I}}(-1)^{\ell (\sigma )} \, x^{L (\sigma )}. 
\end{equation}
We claim that, in our hypotheses, $C_{n}^{I}=C_{n,+}^{I}$.
This is clear if $n \equiv 1 \pmod{2}$ so assume that $n \equiv 0 \pmod{2}$ and $[n-1] \setminus I \subseteq 2 \, {\bf N}$.
Then there exists $\{ a_{1}, \ldots ,a_{s} \} _{<} \subseteq [n-1]$ such that $a_{1} \equiv \ldots \equiv a_{s} \equiv 0 \pmod{2}$
and $I=[1,a_{1}-1]\cup [a_{1}+1,a_{2}-1] \cup \ldots \cup [a_{s}+1, n-1]$. Let $\sigma \in C_{n}^{I}$. Then $\sigma ^{-1}(1)
\in \{ 1,a_{1}+1, \ldots , a_{s}+1 \}$ so $\sigma \in C^{I}_{n,+}$, as desired. $\Box$

Simple examples show that Lemma \ref{chessA} does not hold, in general, if $n \equiv 0 \pmod{2}$ and $[n-1] \setminus I \not \subseteq
2 \,  {\mathbb N}$.

The following simple observation will be used repeatedly in what follows, often without explicit mention.

\begin{lem}
\label{zero}
Let $I\subseteq [n-1]$ and $a\in [2,n-1]$ be such that $[a-2,a+1]\cap I=\emptyset$. Then
\[
\sum_{\scriptsize \begin{array}{l} \{\sigma \in  {S}^{I}_{n}: \\  \sigma(a)=n \} \end{array}}	(-1)^{\ell (\sigma )}x^{L(\sigma )} =	\sum_{\scriptsize \begin{array}{l} \{\sigma \in  {S}^{I}_{n}: \\  \sigma(a)=1 \} \end{array}}	(-1)^{\ell (\sigma )}x^{L(\sigma )}=0. \]
\end{lem}
{\bf Proof:}
In our hypotheses, if $\sigma \in \qa$ then also $\sigma{\check{}}:= \sigma  (a-1,a+1)$ is in the same quotient. Clearly $(\sigma\check{})\check{}=\sigma$ and $|\ell(\sigma)-\ell(\sigma{\check{}})|=1$, while, since $\sigma(a)=n$, $L(\sigma\check{})=L(\sigma)$. Therefore we have that
\begin{align*}\sum_{\scriptsize \begin{array}{l} \{\sigma \in  {S}^{I}_{n}: \\  \sigma(a)=n \} \end{array}} (-1)^{\ell (\sigma )}x^{L(\sigma )}&=& \sum_{\scriptsize \begin{array}{l} \{\sigma \in  {S}^{I}_n: \sigma(a)=n, \\ \sigma (a-1) < \sigma (a+1)  \} \end{array}}\left( (-1)^{\ell (\sigma )}x^{L(\sigma )} + (-1)^{\ell (\sigma\check{} )}x^{L(\sigma\check{} )}\right) =0.\end{align*} 
The proof of the second equality is exactly analogous and is therefore omitted. 
$\Box$

Our next result implies that a connected component of odd cardinality of a subset $I \subseteq [n-1]$ can be shifted to the right, as long as it remains a connected component of $I$, or ``fattened'' by adding the least element bigger than it, without changing the generating function of $(-1)^{\ell(\sigma)} x^{L(\sigma)}$over $\sigma \in S^{I}_{n}$. 

\begin{pro}
\label{scr}
Let $I \subseteq [n-1]$, 
 and $i \in {\mathbb P}$, $k \in {\mathbb N}$ be such that
$[i,i+2k]$ is a connected component of $I$ and $i+2k+2 \not \in I$. 
Then
\begin{equation}\label{sc}\sum_{\sigma \in  S_{n}^{I}}(-1)^{\ell (\sigma )}x^{L(\sigma )} =\sum_{\sigma \in  S_{n}^{I\cup \tilde{I}}}(-1)^{\ell (\sigma )}x^{L(\sigma )}   =\sum_{\sigma \in  S_{n}^{\tilde{I}}}(-1)^{\ell (\sigma )}x^{L(\sigma )}\end{equation}
where $\tilde{I} := (I \setminus \{ i \} ) \cup \{ i+2k+1 \} $.
\end{pro}
{\bf Proof:} Note that
\begin{eqnarray*} \sum_{\sigma \in S_{n}^{I}}(-1)^{\ell (\sigma )}x^{L(\sigma )} & = &\sum_{ \scriptsize  \begin{array}{l}  \{ \sigma \in S_{n}^{I}: \\\sigma (i)>\sigma (i+2k+2)\}  \end{array}}(-1)^{\ell (\sigma )}x^{L(\sigma )} + \sum_{\scriptsize \begin{array}{l} \{ \sigma \in S_{n}^{I}: \; \sigma (i+2k+1)< \\ \; \; \sigma (i+2k+2) \} \end{array}}(-1)^{\ell (\sigma )}x^{L(\sigma )} \\& + & \sum_{j=1}^{2k+1}  \; \; \;  \sum_{\scriptsize \begin{array}{l} \{ \sigma \in S_{n}^{I}: \; \sigma (i+j-1)< \\ \sigma (i+2k+2)<\sigma (i+j)\} \end{array} }(-1)^{\ell (\sigma )}x^{L(\sigma )}.\end{eqnarray*}
Let $j \in [k]$. Then we have that
\[\sum_{\scriptsize \begin{array}{l} \{\sigma \in  {S}^{I}_{n}: \sigma (i+2j-1) < \\  \sigma (i+2k+2) < \sigma (i+2j) \} \end{array}}(-1)^{\ell (\sigma )}x^{L(\sigma )} +\sum_{\scriptsize \begin{array}{l} \{\sigma \in  {S}^{I}_{n}: \sigma (i+2j) < \\  \sigma (i+2k+2) < \sigma (i+2j+1) \} \end{array}}(-1)^{\ell (\sigma )}x^{L(\sigma )} \]
\[
= \sum_{\scriptsize \begin{array}{l} \{\sigma \in  {S}^{I}_{n}: \sigma (i+2j-1) < \\  \sigma (i+2k+2) < \sigma (i+2j) \} \end{array}} [ (-1)^{\ell (\sigma )}x^{L(\sigma )} +(-1)^{\ell (\widetilde{\sigma} )}x^{L(\widetilde{\sigma} )} ]\]where $\tilde{\sigma } := \sigma (i+2j \, , i+2k+2)$.
But $\ell (\widetilde{\sigma})=\ell (\sigma )-1$
and it is easy to see that $L(\widetilde{\sigma}) =L(\sigma )$, so
the above sum is equal to zero. Similarly,
\[\sum_{\scriptsize \begin{array}{l} \{\sigma \in  {S}^{I}_{n}: \sigma (i+2k+2) < \sigma (i) \} \end{array}}(-1)^{\ell (\sigma )}x^{L(\sigma )} +\sum_{\scriptsize \begin{array}{l} \{\sigma \in  {S}^{I}_{n}: \sigma (i) < \\  \sigma (i+2k+2) < \sigma (i+1) \} \end{array}}(-1)^{\ell (\sigma )}x^{L(\sigma )} = 0.\]
Hence
\[ \sum_{\sigma \in S^{I}_{n}}(-1)^{\ell (\sigma )}x^{L(\sigma )} =\sum_{\scriptsize \begin{array}{r} \{\sigma \in  {S}^{I}_{n}: \sigma (i+2k+1) < \sigma (i+2k+2)  \} \end{array}}(-1)^{\ell (\sigma )}x^{L(\sigma )}. \]
This proves the left equality in (\ref{sc}). The proof of the right
equality is exactly analogous and is therefore omitted.
$\Box$

The following is the ``left'' version of Proposition \ref{scr}.

\begin{pro}
\label{scl}
Let $I \subseteq [n-1]$, and $i \in {\mathbb P}$, $k \in {\bf N}$ be such
that $[i+1,i+2k+1]$ is a connected component of $I$ and $i-1 \not \in I$.
Then
\[ \sum_{\sigma \in S^{I}_{n}}(-1)^{\ell (\sigma )}x^{L(\sigma )} = \sum_{\sigma \in S^{I \cup \bar{I}}_{n}}(-1)^{\ell (\sigma )}x^{L(\sigma )}= \sum_{\sigma \in S^{\bar{I}}_{n}}(-1)^{\ell (\sigma )}x^{L(\sigma )} \]
where $\bar{I} := (I \setminus \{ i+2k+1 \} ) \cup \{ i \} $.
\end{pro}
{\bf Proof:} From our hypotheses we have that $[i,i+2k]$ is a connected component of 
$\bar{I}$ and $i+2k+2 \not \in \bar{I}$, so the result follows from 
Proposition \ref{scr}. 
$\Box$

Note that the proofs of the two previous results also prove the following finer 
versions which we also use in the proof of Conjecture \ref{A} in the next section.

\begin{pro}
\label{scrr}
Let $I \subseteq [n-1]$, $i \in {\mathbb P}$, $k \in {\mathbb N}$ be such that
$[i,i+2k]$ is a connected component of $I$, $i+2k+2 \not \in I$, and 
$a \in [n] \setminus (I \cup [i-1,i+2k+2])$ .
Then
\[ \sum_{\scriptsize \begin{array}{r} \{\sigma \in  S_{n}^{I}: \\  \sigma (a)=n \} \end{array}}(-1)^{\ell (\sigma )}x^{L(\sigma )} =\sum_{\scriptsize \begin{array}{r} \{\sigma \in  S_{n}^{I\cup \tilde{I}}: \\  \sigma (a)=n \} \end{array}}(-1)^{\ell (\sigma )}x^{L(\sigma )}   =\sum_{\scriptsize \begin{array}{r} \{\sigma \in  S_{n}^{\tilde{I}}: \\  \sigma (a)=n \} \end{array}}(-1)^{\ell (\sigma )}x^{L(\sigma )}\]where $\tilde{I} := (I \setminus \{ i \} ) \cup \{ i+2k+1 \} $.
\end{pro}

\begin{pro}
\label{sclr}
Let $I \subseteq [n-1]$, $i \in {\mathbb P}$, $k \in {\bf N}$ be such that
$[i+1,i+2k+1]$ is a connected component of $I$, $i-1 \not \in I$, and 
$a \in [n] \setminus (I \cup [i-1,i+2k+2])$ .
Then
\[ \sum_{\scriptsize \begin{array}{r} \{\sigma \in  S_{n}^{I}: \\  \sigma (a)=n \} \end{array}}(-1)^{\ell (\sigma )}x^{L(\sigma )} =\sum_{\scriptsize \begin{array}{r} \{\sigma \in  S_{n}^{I\cup \bar{I}}: \\  \sigma (a)=n \} \end{array}}(-1)^{\ell (\sigma )}x^{L(\sigma )}   =\sum_{\scriptsize \begin{array}{r} \{\sigma \in  S_{n}^{\bar{I}}: \\  \sigma (a)=n \} \end{array}}(-1)^{\ell (\sigma )}x^{L(\sigma )}\]
where $\bar{I} := (I \setminus \{ i+2k+1 \} ) \cup \{ i \} $.
\end{pro}

\section{Main result}

In this section, using the results in the previous one, we obtain closed product formulas for the generating functions of $(-1)^{\ell (\sigma )}
x^{L(\sigma )}$ over the even and odd chessboard elements of any quotient of the symmetric group. In particular, we verify Conjecture \ref{A}.

Let $I \subseteq [n-1]$. We say that $I$ is {\em compressed} if there exists$\{ a_{1}, \ldots, a_{s} \}_{<} \subseteq [n]$ such that $I=[1,a_{1}-1] \cup[a_{1}+1,a_{2}-1] \cup \ldots \cup [a_{s-1}+1,a_{s}-1]$ and $a_{1} \equiv a_{2}\equiv \ldots \equiv a_{s}  \equiv 0 \pmod{2}$. 

\begin{thm}\label{thmA}
Let $n\in {\mathbb P}$, $I \subseteq [n-1]$, and $I_{1},\ldots , I_{s}$ be the connected components of $I$. Then

\begin{equation}
 \sum_{\sigma \in C^{I}_{n,\, +}}(-1)^{\ell (\sigma )}x^{L(\sigma )} ={\left[ \begin{array}{c} \mt \\ \left\lfloor \frac{|I_{1}|+1}{2} \right\rfloor , \ldots , \left\lfloor \frac{|I_{s}|+1}{2} \right\rfloor \end{array} \right]_{x^2}\; \prod_{k= \mt+1}^{\lfloor \frac{n-1}{2} \rfloor} (1-x^{2k}) \,;}\label{c2mp} \end{equation} 
and \begin{equation}
\sum_{\sigma \in C^{I}_{2m,\, -}}(-1)^{\ell (\sigma )}x^{L(\sigma )} =\left\{ \begin{array}{ll}0, & \mbox{if $I$ is compressed and $2m-1 \in I$,} \\ -x^{m} \; {\displaystyle \sum_{\sigma \in C^{I}_{2m,\, +}}}(-1)^{\ell (\sigma )}x^{L(\sigma )}, & \mbox{otherwise,}\end{array} \right.  \label{c2mm} 
\end{equation}
where $\tilde{m} := \sum_{k=1}^{s} \left\lfloor \frac{|I_{k}|+1}{2}\right\rfloor$.
\end{thm}
\vspace{8mm}
{\bf Proof:}
We proceed by induction on $n \in {\mathbb P}$. By
 repeated application of Proposition \ref{scl} we may assume that there exists
 $\{ a_{1},\ldots , a _{s} \} _{<} \subseteq [n]$ such that $I=[1,a_{1}-1]\cup [a_{1}+1,a_{2}-1] \cup \ldots \cup[a_{s-1}+1,a_{s}-1]$, and $|[1,a_{1}-1]| \equiv |[a_{1}+1,a_{2}-1]| \equiv \ldots\equiv |[a_{s-1}+1,a_{s}-1]| \equiv 1 \pmod{2}$, so $a_{1} \equiv a_{2} \equiv \ldots\equiv a_{s} \equiv 0 \pmod{2}$. We have a few cases to distinguish.

\begin{itemize}
\item[i)] Let $n=2m+1$. 
\end{itemize}

We prove \eqref{c2mp}.
\p  If $a_{s}=2m$, then $\sigma(2m+1)=2m+1$ for any $\sigma \in C_{2m+1,+}^{I}$ so\[\sum_{\sigma \in {C_{2m+1, \, +}^{I}}} (-1)^{\ell(\sigma)} x^{L(\sigma)}=\sum_{\tau \in {C_{2m,+}^{I}}} (-1)^{\ell(\bar{\tau})} x^{L(\bar{\tau})}\]
where $\bar{\tau}:= [\tau(1),\cdots,\,\tau(2m),\,2m+1]$. Clearly $\ell(\bar{\tau})=\ell(\tau)$ and $L(\bar{\tau})=L(\tau)$. 

Thus, by our induction hypotheses we conclude that
\[\sum_{\sigma \in {C_{2m+1, \, +}^{I}}} (-1)^{\ell(\sigma)} x^{L(\sigma)}=\sum_{\tau \in {C_{2m,+}^{I}}} (-1)^{\ell(\tau)} x^{L(\tau)}=\mc \] 
as desired.

Assume now that $a_{s}<2m$, that is $a_s\leq 2m-2$. By repeated application of Proposition \ref{scr} we have that

\[ \sum_{\sigma \in C^{I}_{2m+1}}(-1)^{\ell (\sigma )}x^{L(\sigma )} =\sum_{\sigma \in C^{\tilde{I}}_{2m+1}}(-1)^{\ell (\sigma )}x^{L(\sigma )} \]
where $\tilde{I} := [2,a_{1}+1] \cup [a_{1}+3,a_{2}+1] \cup\ldots \cup [a_{s-1}+3,a_{s}+1]$.

Consider first the case $a_{s} < 2m-2$. Then $\sigma^{-1}(2m+1) \in \{ 1, a_s +3 , a_s +5 , \ldots , 2m+1 \}$ for all $\sigma \in C^{\tilde{I}}_{2m+1, \, +}$ so
by Lemma \ref{zero} we have that
\begin{eqnarray*} \sum_{\sigma \in C^{\tilde{I}}_{2m+1}}(-1)^{\ell (\sigma )}x^{L(\sigma )}&=& \sum_{\scriptsize \begin{array}{l} \{\sigma \in C^{\tilde{I}}_{2m+1}: \\\; \; \sigma (1)=2m+1 \} \end{array}}(-1)^{\ell (\sigma )}x^{L(\sigma )}+ \sum_{\scriptsize \begin{array}{l} \{\sigma \in C^{\tilde{I}}_{2m+1}: \\\; \; \sigma (a_{s}+3)=2m+1 \} \end{array}}(-1)^{\ell (\sigma )}x^{L(\sigma )} \\ & &+ \sum_{\scriptsize \begin{array}{l} \{\sigma \in C^{\tilde{I}}_{2m+1}: \\\; \; \sigma (2m+1)=2m+1 \} \end{array}}(-1)^{\ell (\sigma )}x^{L(\sigma )}. \end{eqnarray*}
But 
\begin{eqnarray*} \sum_{\scriptsize \begin{array}{l} \{\sigma \in C^{\tilde{I}}_{2m+1}: \\\; \; \sigma (2m+1)=2m+1 \} \end{array}}(-1)^{\ell (\sigma )}x^{L(\sigma )} & =& \sum_{\tau \in C^{\tilde{I}}_{2m,+}}(-1)^{\ell (\tau )}x^{L(\tau )} \\& = & \prod_{k=\mt +1}^{m-1} (1-x^{2k}) \mc,\end{eqnarray*}
by our induction hypotheses, while

\[ \sum_{\scriptsize \begin{array}{l} \{\sigma \in C^{\tilde{I}}_{2m+1}: \\\; \; \sigma (1)=2m+1 \} \end{array}}(-1)^{\ell (\sigma )}x^{L(\sigma )} = \sum_{\tau \in C^{\bar{I}}_{2m,-}}(-1)^{\ell (\bar{\tau} )}x^{L(\bar{\tau} )} ,\]
where $\bar{\tau}:= [2m+1,\tau(1),\cdots,\,\tau(2m)]$ and $\bar{I}=[1,a_{1}]\cup[a_1 + 2, a_2]\cup\cdots\cup [a_{s-1}+2,a_s]$.
But $\ell(\bar{\tau})=\ell(\tau)+2m$ and $L(\bar{\tau})=L(\tau)+m$, hence, by our induction hypotheses:
\begin{eqnarray*} \sum_{\scriptsize \begin{array}{l} \{\sigma \in C^{\tilde{I}}_{2m+1}: \\\; \; \sigma (1)=2m+1 \} \end{array}}(-1)^{\ell (\sigma )}x^{L(\sigma )} & =& x^m \sum_{\tau \in C^{\bar{I}}_{2m,-}}(-1)^{\ell (\tau )}x^{L(\tau)} =\\&=& -x^{2m} \prod_{k=\mt +1}^{m-1} (1-x^{2k}) \mc\end{eqnarray*}

(note that $\frac{a_s}{2}=\mt$). Finally, by repeated application of Proposition 
\ref{sclr} we have 
\[  \sum_{\scriptsize \begin{array}{l} \{\sigma \in C^{\tilde{I}}_{2m+1}: \\\; \; \sigma (a_{s}+3)=2m+1 \} \end{array}}(-1)^{\ell (\sigma )}x^{L(\sigma )} = \sum_{\scriptsize \begin{array}{l} \{\sigma \in C^{I}_{2m+1}: \\\; \; \sigma (a_{s}+3)=2m+1 \} \end{array}}(-1)^{\ell (\sigma )}x^{L(\sigma )} =0 ,\]
by Lemma \ref{zero}. So
\begin{eqnarray*}\sum_{\sigma \in C^{I}_{2m+1}}(-1)^{\ell (\sigma )}x^{L(\sigma )} &=&\prod_{k=\mt+1}^{m} (1-x^{2k}) \mc,\end{eqnarray*}
as desired.

If $a_{s}=2m-2$ then we similarly have that
\begin{eqnarray*}\sum_{\sigma \in C^{\tilde{I}}_{2m+1}}(-1)^{\ell (\sigma )}x^{L(\sigma )} &=&  \sum_{\scriptsize \begin{array}{l} \{\sigma \in C^{\tilde{I}}_{2m+1}: \\\; \; \sigma (2m+1)=2m+1 \} \end{array}}(-1)^{\ell (\sigma )}x^{L(\sigma )} + \\& +& \sum_{\scriptsize \begin{array}{l} \{\sigma \in C^{\tilde{I}}_{2m+1}: \\\; \; \sigma (1)=2m+1 \} \end{array}}(-1)^{\ell (\sigma )}x^{L(\sigma )},\end{eqnarray*}
and the result follows exactly as before. This proves (\ref{c2mp}).

\begin{itemize}
\item[ii)] Let $n=2m$.
\end{itemize}

We first prove \eqref{c2mp}. Assume first that $a_{s}=2m$.
Then $\sigma^{-1}(2m) \in \{ a_1 , \ldots , a_s \}$ for any $\sigma \in C^{I}_{2m, \, +}$ so 
\begin{equation}\label{comp1} \sum_{\sigma \in C^{I}_{2m,\, +}}(-1)^{\ell (\sigma )}x^{L(\sigma )}  = \sum_{j=1}^{s} \; \sum_{\scriptsize\begin{array}{l} \{\sigma \in C^{I}_{2m,\, +}: \\\; \; \sigma (a_{j})=2m \} \end{array}}(-1)^{\ell (\sigma )}x^{L(\sigma )} .\end{equation}
Let $j \in [s]$. Then $\sigma^{-1}(2m-1)=a_j -1$ for all $\sigma \in C^{I}_{2m, \, +}$ such that $\sigma(a_j)=2m$ so
\begin{eqnarray*}\label{comp2} \sum_{\scriptsize \begin{array}{l} \{\sigma \in C^{I}_{2m,\, +}: \\\; \; \sigma (a_{j})=2m \} \end{array}}(-1)^{\ell (\sigma )}x^{L(\sigma )} &=&  \sum_{\scriptsize \begin{array}{l} \{\sigma \in C^{I}_{2m,\, +}: \\\; \; \sigma (a_{j})=2m \\\; \; \sigma(a_{j}-1)=2m-1 \} \end{array}}(-1)^{\ell (\sigma )}x^{L(\sigma )} \\& = & \sum_{\tau \in C^{\tilde{I}_{j}}_{2m-2,\, +}}(-1)^{\ell ({\bar{\tau}})}x^{L({\bar{\tau}} )} ,\end{eqnarray*}where $\bar{\tau} := [\tau (1), \ldots , \tau (a_{j}-2), \,2m-1, \; 2m, \; \tau(a_{j}-1), \ldots , \tau (2m-2)]$ and $ \tilde{I}_{j} := [1,a_{1}-1] \cup [a_{1}+1,a_{2}-1] \cup \cdots \cup [a_{j-2}+1,a_{j-1}-1] \cup [a_{j-1}+1,a_{j}-3]\cup [a_{j}-1,a_{j+1}-3]  \cup  [a_{j+1}-1,a_{j+2}-3]\cup\cdots \cup [a_{s-1}-1,2m-3]$. But $\ell (\bar{\tau})=\ell (\tau )+2(2m-a_{j})$ and $L(\bar{\tau})=L(\tau )+(2m-a_{j})$ so we conclude from (\ref{comp1}), (\ref{comp2}), and our induction hypotheses that
\begin{eqnarray*}\sum_{\sigma \in C^{I}_{2m,\, +}}(-1)^{\ell (\sigma )}x^{L(\sigma )} & = & \sum_{j=1}^{s} x^{2m-a_{j}} \; \sum_{\tau \in C^{\tilde{I}_{j}}_{2m-2,\, +}}(-1)^{\ell (\tau )}x^{L(\tau )} \\& =& \sum_{j=1}^{s} x^{2m-a_{j}} \left[ \scriptsize{ \begin{array}{c}m-1 \\\lfloor \frac{|I_1|+1}{2} \rfloor, \ldots , \lfloor \frac{|I_{j-1}|+1}{2} \rfloor,\lfloor \frac{|I_j|-1}{2} \rfloor,\lfloor \frac{|I_{j+1}|+1}{2} \rfloor,\ldots , \lfloor \frac{|I_s|+1}{2} \rfloor
\end{array}} \right]_{x^{2}} \\& = & \left[ \scriptsize{ \begin{array}{c}m \\\lfloor \frac{|I_1|+1}{2} \rfloor, \ldots , \lfloor \frac{|I_s|+1}{2} \rfloor\end{array}} \right]_{x^{2}} ,\end{eqnarray*}
and the result again follows.

Assume now that $a_{s}<2m$. Then $a_{s} \leq 2m-2$ and by repeated application of Proposition \ref{scr} we have that
\[ \sum_{\sigma \in C^{I}_{2m,\, +}}(-1)^{\ell (\sigma )}x^{L(\sigma )} =\sum_{\sigma \in C^{\tilde{I}}_{2m,\, +}}(-1)^{\ell (\sigma )}x^{L(\sigma )} \]where $\tilde{I} := [1,a_{1}] \cup [a_{1}+2,a_{2}] \cup\ldots \cup [a_{s-1}+2,a_{s}]$.

If $a_{s} < 2m-2$, then $\sigma^{-1}(2m) \in \{ a_s +2 , a_s +4 , \ldots , 2m \}$ 
for all $\sigma \in C^{\tilde{I}}_{2m,\, +}$ so by Lemma \ref{zero} have that
\[ \sum_{\sigma \in C^{\tilde{I}}_{2m,\, +}}(-1)^{\ell (\sigma )}x^{L(\sigma )}= \sum_{\scriptsize \begin{array}{l} \{\sigma \in C^{\tilde{I}}_{2m,\, +}: \\\; \; \sigma (a_{s}+2)=2m \} \end{array}}(-1)^{\ell (\sigma )}x^{L(\sigma )} + \sum_{\scriptsize \begin{array}{l} \{\sigma \in C^{\tilde{I}}_{2m,\, +}: \\\; \; \sigma (2m)=2m \} \end{array}}(-1)^{\ell (\sigma )}x^{L(\sigma )}. \]
Now, by our induction hypotheses,
\begin{eqnarray*} \sum_{\scriptsize \begin{array}{l} \{\sigma \in C^{\tilde{I}}_{2m,\, +}: \\\; \; \sigma (2m)=2m \} \end{array}}(-1)^{\ell (\sigma )}x^{L(\sigma )} & =& \sum_{\tau \in C^{\tilde{I}}_{2m-1}}(-1)^{\ell (\tau )}x^{L(\tau )} \\& = & \prod_{k=\frac{a_{s}+2}{2}}^{m-1} (1-x^{2k})\left[  \begin{array}{c}\frac{a_{s}}{2} \\\lfloor \frac{|I_1|+1}{2} \rfloor , \ldots , \lfloor \frac{|I_s|+1}{2} \rfloor\end{array} \right]_{x^{2}} .\end{eqnarray*}
Also, by repeated application of Proposition \ref{sclr} we get
\[  \sum_{\scriptsize \begin{array}{l} \{\sigma \in C^{\tilde{I}}_{2m,\, +}: \\\; \; \sigma (a_{s}+2)=2m \} \end{array}}(-1)^{\ell (\sigma )}x^{L(\sigma )} = \sum_{\scriptsize \begin{array}{l} \{\sigma \in C^{I}_{2m,\, +}: \\\; \; \sigma (a_{s}+2)=2m \} \end{array}}(-1)^{\ell (\sigma )}x^{L(\sigma )} =0 \]
by Lemma \ref{zero} and the result again follows.

If $a_{s}=2m-2$ then we have similarly that
\begin{eqnarray*}\sum_{\sigma \in C^{\tilde{I}}_{2m,\, +}}(-1)^{\ell (\sigma )}x^{L(\sigma )} = \sum_{\scriptsize \begin{array}{l} \{\sigma \in C^{\tilde{I}}_{2m,\, +}: \\\; \; \sigma (2m)=2m \} \end{array}}(-1)^{\ell (\sigma )}x^{L(\sigma )}  = \sum_{\tau \in C^{\tilde{I}}_{2m-1}}(-1)^{\ell (\tau )}x^{L(\tau )} ,\end{eqnarray*}
and the result follows exactly as above. This proves (\ref{c2mp}).

\vspace{5mm}
\p
We now prove (\ref{c2mm}). 

If $a_{s}=2m$ then $C^{I}_{2m, \, -}=\emptyset$ so (\ref{c2mm}) clearly
holds. So assume that $a_{s} <2m$. Then $a_{s} \leq 2m -2$ and $\sigma^{-1}(1) \in \{ a_s +2 , a_s +4 , \ldots , 2m \}$ 
for all $\sigma \in C^{I}_{2m,\, -}$ so by Lemma \ref{zero} we have that
\begin{eqnarray*}\sum_{\sigma \in C^{I}_{2m,\, -}}(-1)^{\ell (\sigma )}x^{L(\sigma )} = \sum_{\scriptsize \begin{array}{l} \{\sigma \in C^{I}_{2m,\, -}: \\\; \; \sigma (2m)=1 \} \end{array}}(-1)^{\ell (\sigma )}x^{L(\sigma )}  =\sum_{\tau \in C^{I}_{2m-1}}(-1)^{\ell (\check{\tau })}x^{L(\check{\tau })} ,\end{eqnarray*}
where  \p $\check{\tau } := [ \tau (1)+1, \; \tau (2)+1, \ldots ,\tau (2m-1) +1, \; 1]$. But, $\ell (\check{\tau })=\ell (\tau )+2m-1$, and$L(\check{\tau })=L(\tau )+m$, so \begin{eqnarray*}\sum_{\tau \in C^{I}_{2m-1}}(-1)^{\ell (\check{\tau })}x^{L(\check{\tau })} & = &-x^{m} \sum_{\tau \in C^{I}_{2m-1}}(-1)^{\ell (\tau )}x^{L(\tau )} \\& = & -x^{m} \; \left[  \begin{array}{c}\frac{a_{s}}{2} \\\lfloor \frac{|I_1|+1}{2} \rfloor , \ldots , \lfloor \frac{|I_s|+1}{2} \rfloor\end{array} \right]_{x^{2}} \prod_{k = \frac{a_{s}+2}{2}}^{m-1} (1-x^{2k}) ,\end{eqnarray*}
by our induction hypotheses, and the result follows from (\ref{c2mp}). 

This concludes the induction step and hence the proof.
$\Box$

As a corollary of Theorem \ref{thmA} we obtain a proof of Conjecture \ref{A} (i.e., of Conjecture C of \cite{KV}).
\begin{thm}
Let $n \in {\mathbb P}$, $I \subseteq [n-1]$, and $I_{1}, \ldots , I_{s}$ be the connected components of $I$. Then
\begin{align}\small\sum_{\sigma \in C^{I}_{2m+1}} (-1)^{\ell (\sigma )} \chi(\sigma) x^{L(\sigma )} = &\mc \displaystyle\prod\limits_{\tiny{k=\mt + 1}}^{\tiny m}{(1-x^{2k})} \\\sum_{\sigma \in C^{I}_{2m}} (-1)^{\ell (\sigma )} \chi(\sigma) x^{L(\sigma )} = &\left\{ \begin{array}{ll}\mc , & \mbox{if $m=\mt$,} \\ (1+x^{m}) \mc \displaystyle\prod\limits_{\tiny{k=\mt + 1}}^{\tiny {m-1}}{(1-x^{2k})}, & \mbox{otherwise,}\end{array} \right.  
\end{align}
\normalsize where $\tilde{m} := \sum_{k=1}^{s} \left\lfloor \frac{|I_{k}|+1}{2} \right\rfloor$.
\end{thm}
{\bf Proof:} The first equation follows immediately from (\ref{c2mp}) of Theorem \ref{thmA} since $C^{I}_{2m+1}=C^{I}_{2m+1,+}$ and $\chi$ is trivial on $C^{I}_{2m+1,+}$.
Also, by definition of $\chi$,
\[ \sum_{\sigma \in C_{2m}^{I}} (-1)^{\ell (\sigma )}\chi (\sigma ) x ^{L (\sigma )} = \sum_{\sigma \in C_{2m,+}^{I}} (-1)^{\ell (\sigma )}
x^{L (\sigma )}
- \sum_{\sigma \in C_{2m , -}^{I}} (-1)^{\ell (\sigma )} x^{L(\sigma )} \]
so the second equation also follows immediately from Theorem \ref{thmA} and the observation that $m=\tilde{m}$ if and only if $I$ is compressed and $2m -1 \in I$.
$\Box$

Also as an immediate consequence of Theorem \ref{thmA} we obtain closed product formulas for the generating function of $(-1)^{\ell (\sigma )}
x^{L(\sigma )}$ over any quotient of $S_{n}$.
\begin{cor}
Let $n \in {\mathbb P}$, $I \subseteq [n-1]$, and $I_{1}, \ldots , I_{s}$ be the connected components of $I$. Then
 \begin{align}
 \sum_{\sigma \in S_{n}^{I}} (-1)^{\ell (\sigma )} x^{L(\sigma )} &=\left\{ \begin{array}{l}\mc 
{\displaystyle \prod_{k=\tilde{m}+1}^{\left\lfloor \frac{n-1}{2} \right\rfloor }}(1-x^{2k}), \\
 \mbox{if $n \equiv 1 \pmod{2}$, or if $n = 2 \tilde{m}$,} \\ 
 (1+x^{m}) \mc {\displaystyle \prod_{k=\tilde{m}+1}^{\left\lfloor \frac{n-1}{2} \right\rfloor }}(1-x^{2k}), \\
  \mbox{otherwise,}\end{array} \right.  
 \end{align}
where $\tilde{m} := \sum_{k=1}^{s} \left\lfloor \frac{|I_{k}|+1}{2} \right\rfloor $.
\end{cor}
{\bf Proof:} This follows immediately from Theorem \ref{thmA}, the definition of $C_{n}$, and the fact that equation (\ref{chessA0}) holds for all
$n \in \mathbb P$ and $I \subseteq [n-1]$. $\Box$

In particular, we obtain the following result for the whole group. 
\begin{cor}
Let $n \in {\mathbb P}$. Then
\[\sum_{\sigma \in S_n} (-1)^{\ell(\sigma)}x^{L(\sigma)} = \left\{ \begin{array}{lll}  &\displaystyle\prod\limits_{\tiny{j=1}}^{\tiny m}(1-x^{2j}), & \mbox{ if } n=2m+1, \\ \\& (1-x^m) \displaystyle\prod\limits_{\tiny{j=1}}^{\tiny m-1}(1-x^{2j}),  & \mbox{ if } n=2m.\end{array}\right.\]
\end{cor}

\section{Type B quotients}

In this section, using Theorem \ref{thmA},
we prove Conjecture \ref{B}. A different proof of this conjecturte
appears in \cite{Lan} (see also \cite{Car}).

Our first result is the analogue, for the odd length function $L$,
of a well known description of the ordinary length function of the hyperoctahedral group (see, e.g., \cite[(8.1)]{BB}). Its proof is a simple verification and is omitted.

Given $\sigma \in B_n$ we let
\begin{align*}\p oinv(\sigma):=& |\{(i,j)\in [n]\times[n] \; : \; i<j,\,\sigma(i)>\sigma(j),\,i\not\equiv j \pmod{2}\}|, \\oneg(\sigma):=& |\{ i\in [n]\; : \; \sigma(i)<0,\,i\not\equiv 0 \pmod{2}\}|, \\onsp(\sigma):=& |\{(i,j)\in [n]\times[n] \; : \; \sigma(i)+\sigma(j)<0,\,i\not\equiv j \pmod{2}\}|.\end{align*}
\begin{pro}\label{LB}
Let $\sigma \in B_n$. Then 
\[
L(\sigma)=oinv(\sigma)+oneg(\sigma)+onsp(\sigma).\]
\end{pro}   
  
Note that the previous result is similar to, but different from, Lemma 6 of \cite{VS}.

The next result is the analogue, for type $B$, of Proposition \ref{scr}. Its proof is identical, ``mutatis mutandis'',
to that of Proposition \ref{scr} and is therefore omitted.
\begin{pro}
\label{shiftB}
Let $I \subseteq [0,n-1]$, and $i \in {\mathbb P}$, $k \in {\mathbb N}$ be such that
$[i,i+2k]$ is a connected component of $I$ and $i+2k+2 \not \in I$. Then
\begin{equation}\label{shiftB1}\sum_{\sigma \in B_{n}^{I}}(-1)^{\ell (\sigma )}x^{L(\sigma )} =\sum_{\sigma \in B_{n}^{I\cup \tilde{I}}}(-1)^{\ell (\sigma )}x^{L(\sigma )} =\sum_{\sigma \in B_{n}^{\tilde{I}}}(-1)^{\ell (\sigma )}x^{L(\sigma )} , \end{equation}
where $\tilde{I} :=  (I \setminus \{ i \} ) \cup \{ i+2k+1 \} $.
\end{pro}
Our next result describes the effect, on  the generating function of $(-1)^{\ell (\sigma )} x^{L(\sigma )}$ over
$B_{n}^{J}$, of ``compressing'' the connected component of $J$ that contains $0$.
\begin{pro}
\label{inflateB}
Let $J\subseteq [0,n-1]$ and $a\in [0,n-1]$ be such that $[0,a-1]\subseteq J$, $a,a+1\,\notin J$. Then
\[\sumb=(1-x^{a+1}) \sum_{\sigma \in {B_n^{J\cup\{a\}}}} (-1)^{\ell(\sigma)} x^{L(\sigma)}.\]
\end{pro}
{\bf Proof:}
Note first that $\sigma(a) \geq 0$ for all $\sigma \in B^{J}_{n}$. Hence we have that
\begin{eqnarray*}\sumb &=& \sum_{\scriptsize \begin{array}{l} \{\sigma \in \qb : \\\;\; \sigma(a+1)>\sigma(a)\}\end{array}} (-1)^{\ell(\sigma)} x^{L(\sigma)} + \sum_{\scriptsize \begin{array}{l} \{\sigma \in {\qb} : \\\;\; \sigma(a+1)<\sigma(-a)\}\end{array}} (-1)^{\ell(\sigma)} x^{L(\sigma)} + \\&+&\sum_{j=1}^{a}\left(\sum_{\scriptsize \begin{array}{l} \{\sigma \in {\qb} : \sigma(j-1)< \\\;\; \sigma(a+1)<\sigma(j)\}\end{array}} (-1)^{\ell(\sigma)} x^{L(\sigma)} +\sum_{\scriptsize \begin{array}{l} \{\sigma \in \qb : \sigma(-j)<\\\;\; \sigma(a+1)<\sigma(-j+1)\}\end{array}} (-1)^{\ell(\sigma)} x^{L(\sigma)} \right).\end{eqnarray*}

We claim that
\begin{equation}
\label{6.1.2}
\sumb = \sum_{\scriptsize \begin{array}{l} \{\sigma \in \qb : \\\;\; \sigma(a+1)>\sigma(a)\}\end{array}} (-1)^{\ell(\sigma)} x^{L(\sigma)} + \sum_{\scriptsize \begin{array}{l} \{\sigma \in {\qb} : \\\;\; \sigma(a+1)<\sigma(-a)\}\end{array}} (-1)^{\ell(\sigma)} x^{L(\sigma)}.
\end{equation}

To show this we have two cases to distinguish.
\begin{itemize}
\item[i)] $a\equiv 1 \pmod{2}$
\end{itemize}
Let $j\in [\frac{a-1}{2}]$. Then we have that 
\[\sum_{\scriptsize \begin{array}{l} \{\sigma \in {\qb} : \sigma(2j-1)<\\\;\; \sigma(a+1)<\sigma(2j)\}\end{array}} (-1)^{\ell(\sigma)} x^{L(\sigma)} +\sum_{\scriptsize \begin{array}{l} \{\sigma \in \qb :  \sigma(2j)<\\\;\;\sigma(a+1)<\sigma(2j+1)\}\end{array}} (-1)^{\ell(\sigma)} x^{L(\sigma)} \]
\begin{eqnarray}\label{5.2.1}&=& \sum_{\scriptsize \begin{array}{l} \{\sigma \in {\qb} : \sigma(2j-1)< \\\;\; \sigma(a+1)<\sigma(2j)\}\end{array}} \left( (-1)^{\ell(\sigma)} x^{L(\sigma)} + (-1)^{\ell(\bar{\sigma})} x^{L(\bar{\sigma})} \right)\end{eqnarray}
where $\bar{\sigma} := \sigma \, (a+1,2j)(-2j \, ,-a-1)$.
But $\ell(\bar{\sigma})=\ell(\sigma)-1$ and, since $a\equiv 1 \pmod{2}$, $L(\bar{\sigma})=L(\sigma)$ so the sum in (\ref{5.2.1}) is equal to $0$. 

Similarly
\[\sum_{\scriptsize \begin{array}{l} \{\sigma \in {\qb} : \sigma(-2j-1)< \\\;\; \sigma(a+1)<\sigma(-2j)\}\end{array}} (-1)^{\ell(\sigma)} x^{L(\sigma)} +\sum_{\scriptsize \begin{array}{l} \{\sigma \in \qb : \sigma(-2j)< \\\;\; \sigma(a+1)<\sigma(-2j+1)\}\end{array}} (-1)^{\ell(\sigma)} x^{L(\sigma)} \]
\begin{eqnarray}\label{5.2.2}&=& \sum_{\scriptsize \begin{array}{l} \{\sigma \in {\qb} : \sigma(-2j-1)<\\\;\; \sigma(a+1)<\sigma(-2j)\}\end{array}} ((-1)^{\ell(\sigma)} x^{L(\sigma)} + (-1)^{\ell(\bar{\sigma})} x^{L(\bar{\sigma})})\end{eqnarray}
where $\bar{\sigma} := \sigma (-2j,a+1)(-a-1,2j)$.
Again, $\ell(\bar{\sigma})=\ell(\sigma)-1$ and $L(\bar{\sigma})=L(\sigma)$ so the sum in (\ref{5.2.2}) is equal to $0$. 

Furthermore
\[
\sum_{\scriptsize \begin{array}{l} \{\sigma \in \qb : \\\;\; 0<\sigma(a+1)<\sigma(1)\}\end{array}} (-1)^{\ell(\sigma)} x^{L(\sigma)} + \sum_{\scriptsize \begin{array}{l} \{\sigma \in {\qb} : \\\;\; -\sigma(1)<\sigma(a+1)<0\}\end{array}} (-1)^{\ell(\sigma)} x^{L(\sigma)} \]
\begin{eqnarray}\label{5.2.3}&=& \sum_{\scriptsize \begin{array}{l} \{\sigma \in \qb : \\\;\; 0<\sigma(a+1)<\sigma(1)\}\end{array}} ((-1)^{\ell(\sigma)} x^{L(\sigma)} + (-1)^{\ell(\bar{\sigma})} x^{L(\bar{\sigma})})\end{eqnarray}
where $\bar{\sigma} := \sigma (a+1,-a-1)$. 
Clearly $\ell(\bar{\sigma})=\ell(\sigma)+1$, while, since $a\equiv 1 \pmod{2}$, $L(\bar{\sigma})=L(\sigma)$, so the sum in (\ref{5.2.3}) is also equal to $0$.

\begin{itemize}
\item[ii)] $a\equiv 0 \pmod{2}$
\end{itemize}
If $a=0$ then (\ref{6.1.2}) is clear, so assume $a \geq 1$.
Let $j\in [\frac{a}{2}]$. Then we similarly have that
\[\sum_{\scriptsize \begin{array}{l} \{\sigma \in {\qb} : \sigma(2j-2)<\\\;\; \sigma(a+1)<\sigma(2j-1)\}\end{array}} (-1)^{\ell(\sigma)} x^{L(\sigma)} +\sum_{\scriptsize \begin{array}{l} \{\sigma \in \qb : \sigma(2j-1)<\\\;\; \sigma(a+1)<\sigma(2j)\}\end{array}} (-1)^{\ell(\sigma)} x^{L(\sigma)} \]
\begin{eqnarray}\label{5.3.7}&=& \sum_{\scriptsize \begin{array}{l} \{\sigma \in {\qb} : \sigma(2j-2)<\\\;\; \sigma(a+1)<\sigma(2j-1)\}\end{array}} {\left((-1)^{\ell(\sigma)} x^{L(\sigma)} +(-1)^{\ell(\bar{\sigma})} x^{L(\bar{\sigma})}\right)}\end{eqnarray}
where 
$\bar{\sigma} := \sigma (2j-1,a+1)(-2j+1,-a-1)$. 
Since $\ell(\bar{\sigma})=\ell(\sigma)-1$ and, since $a\equiv 0 \pmod{2}$, $L(\bar{\sigma})=L(\sigma)$ the sum in (\ref{5.3.7}) is equal to $0$. 

Similarly
\[\sum_{\scriptsize \begin{array}{l} \{\sigma \in {\qb} :\sigma(-2j)< \\\;\; \sigma(a+1)<\sigma(-2j+1)\}\end{array}} (-1)^{\ell(\sigma)} x^{L(\sigma)} +\sum_{\scriptsize \begin{array}{l} \{\sigma \in \qb : \sigma(-2j+1)<\\\;\; \sigma(a+1)<\sigma(-2j+2)\}\end{array}} (-1)^{\ell(\sigma)} x^{L(\sigma)} \]
\begin{eqnarray}\sum_{\scriptsize \begin{array}{l} \{\sigma \in \qb : \sigma(-2j)<\\\;\; \sigma(a+1)<\sigma(-2j+1)\}\end{array}}{\left((-1)^{\ell(\sigma)} x^{L(\sigma)}+(-1)^{\ell(\bar{\sigma})} x^{L(\bar{\sigma})}\right)}= 0\end{eqnarray}
where 
$\bar{\sigma} := \sigma (a+1,-2j+1)(-a-1,2j-1)$.

This proves our claim. Therefore we have from (\ref{6.1.2}) that
\begin{eqnarray*}\sumb&=& \sum_{\scriptsize \begin{array}{l} \{\sigma \in \qb : \\\;\; \sigma(a+1)>\sigma(a)\}\end{array}} (-1)^{\ell(\sigma)} x^{L(\sigma)} + \sum_{\scriptsize \begin{array}{l} \{\sigma \in {\qb} : \\\;\; \sigma(a+1)<\sigma(-a)\}\end{array}} (-1)^{\ell(\sigma)} x^{L(\sigma)} \\&=&\sum_{\scriptsize \begin{array}{l} \{\sigma \in \qb : \\\;\; \sigma(a)<\sigma(a+1)\}\end{array}} {\left((-1)^{\ell(\sigma)} x^{L(\sigma)} +(-1)^{\ell(\bar{\sigma})} x^{L(\bar{\sigma})}\right)}\end{eqnarray*}
where $\bar{\sigma} := \sigma (-a-1,a+1)$.
But $\ell(\bar{\sigma})=\ell(\sigma)+2a+1$ and, by Proposition \ref{LB}, $L(\bar{\sigma})=L(\sigma)+a+1$ (note that this is true for two different reasons
depending on the parity of $a$),  
therefore
\[
\sumb=(1-x^{a+1})\sum_{\scriptsize \begin{array}{l} \{\sigma \in \qb : \\\;\; \sigma(a)<\sigma(a+1)\}\end{array}} (-1)^{\ell(\sigma)} x^{L(\sigma)}
\]
and the result follows. 
$\Box$

We can now prove Conjecture \ref{B}. 
For $J\subseteq [0,n-1]$ we define $J_0\subseteq J$ to be the connected component of $J$ which contains $0$, if $0\in J$, or $J_{0} := \emptyset$ otherwise. Let $J_1,\ldots,J_s$ be the remaining ordered connected components.

\begin{thm}
\label{typeB}
Let $n\in {\mathbb P}$, $J \subseteq [0,n-1]$, and $J_0,\ldots, J_s$ be the connected
components of $J$ indexed as just described.  Then
\begin{equation}\sumb=\frac{\prod_{j=a+1}^{n}(1-x^j)}{\prod_{i=1}^{\mt}(1-x^{2i})}  \mcb \end{equation}
where $\mt:=\sum_{i=1}^s \left\lfloor \frac{|J_{i}|+1}{2}\right\rfloor$ and 
$a:=\min\{ [0, \,n-1] \setminus J\}$.
\end{thm}
{\bf Proof:}
We distinguish the cases $n$ even and $n$ odd.

Let $n=2m+1$ and suppose first that $J_{0}=\emptyset$. Then from Propositions \ref{fact2} and \ref{ascending}, and Theorem \ref{thmA},
we have that 
\begin{eqnarray*} \sumb&=&\left(\sum_{\sigma \in {B_n^{[n-1]}}}{(-1)^{\ell(\sigma)} x^{L(\sigma)}}\right)\left(\sum_{\sigma \in {S_n^{J}}}{(-1)^{\ell(\sigma)} x^{L(\sigma)}}\right) \\ &=&\prod_{j=1}^{m+1}{(1-x^{2j-1})} \mcb \prod_{i=\mt+1}^{m}{(1-x^{2i})} \end{eqnarray*}
and the result follows.
Suppose now that $0\in J,$ say $J_0=[0,\,a-1]$. Then by repeated application of Proposition \ref{inflateB} we have that
\begin{equation}\sumb=\frac{1}{\prod_{i=1}^{a} (1-x^i)}\sum_{\sigma \in {B_n^{J\setminus J_0}}} (-1)^{\ell(\sigma)} x^{L(\sigma)},\end{equation}
and the result follows from the previous case.

Let now $n=2m$. By repeated application of Proposition \ref{shiftB} we may assume that there
exist $\{ a_1, \ldots , a_{s} \}_{<} \subseteq [0,n-2]$ such that 
\[ J_1=[a_1 +1,\,a_2-1],\;J_2=[a_2 + 1,\,a_3-1],\ldots ,J_s=[a_s+1,\,n-1]\]
and $ a_1 \equiv \cdots \equiv a_s \equiv 0 \pmod 2$.
Let $\widetilde{a}:= m-\mt=a_1 /2$, $\widetilde{J}_{0}:= [0, a_{1}-1]$, 
and $\widetilde{J}:= \widetilde{J}_{0} \cup J_{1} \cup \ldots \cup J_{s}$. 
Then, by Propositions \ref{ascending}, \ref{fact2}, and Theorem \ref{thmA}
\begin{eqnarray*} \sum_{\sigma \in B^{\widetilde{J}\setminus\{0\}}_{n}}{(-1)^{\ell({\sigma})}x^{L(\sigma)}}&=&\left(\sum_{\sigma \in {B_n^{[n-1]}}}{(-1)^{\ell(\sigma)} x^{L(\sigma)}}\right)\left(\sum_{\sigma \in S^{\widetilde{J}\setminus \{0\}}_n} {(-1)^{\ell(\sigma)} x^{L(\sigma)}}\right) \\
&=& \prod_{j=1}^{m} (1-x^{2j-1}) {\left[ \begin{array}{c} m \\\widetilde{a},\,\left\lfloor \frac{|J_{1}|+1}{2} \right\rfloor , \ldots , \left\lfloor \frac{|J_{s}|+1}{2} \right\rfloor \end{array} \right]_{x^2}} \\
&=& \prod_{j=1}^{m} (1-x^{2j-1}) \frac{[m]_{\tiny x^2}!}{[\widetilde{a}]_{\tiny x^2}!\,[\mt]_{\tiny x^2}!} \mcb .\end{eqnarray*}
But, by repeated application of Proposition \ref{inflateB}
\[
\sumb=\prod_{i=a+1}^{a_1}(1-x^i) \sum_{\sigma \in B^{\widetilde{J}}_n} (-1)^{\ell(\sigma)}x^{L(\sigma)}\]
and, by Proposition \ref{2.8b} 
\[
\sum_{\sigma \in B^{\widetilde{J}\setminus\{0\}}_{n}}{(-1)^{\ell({\sigma})}x^{L(\sigma)}}=\left(\sum_{\sigma \in B^{\widetilde{J}}_{n}}{(-1)^{\ell({\sigma})}x^{L(\sigma)}}\right)\left(\sum_{\sigma\in B^{[a_{1}-1]}_{a_{1}}}{(-1)^{\ell({\sigma})}x^{L(\sigma)}}\right).\]
Combining the previous identities we get
\begin{eqnarray*}\sumb&=&\frac{\prod_{j=1}^{m}(1-x^{2j-1}) \prod_{i=a+1}^{2\widetilde{a}}(1-x^i)}{\prod_{i=1}^{\widetilde{a}}(1-x^{2i-1})} \frac{[m]_{\tiny x^2}!}{[\widetilde{a}]_{\tiny x^2}!\,[\mt]_{\tiny x^2}!} \mcb \\&=&\frac{\prod_{j=a+1}^{n}(1-x^i)}{\prod_{i=1}^{\mt}(1-x^{2i})}  \mcb ,\end{eqnarray*}
as desired. $\Box$

\end{document}